\documentclass[11pt]{amsart}
\usepackage{amsmath,amssymb,amsfonts,amsthm}
\usepackage[left=1 in, right=1 in, top=1.5 in, bottom=1.5 in]{geometry}
\usepackage{graphicx} 
\usepackage{cleveref}
\usepackage{units}
\usepackage{ytableau}
\usepackage{bm}
\usepackage{mathtools}

\title{Super major index and Thrall's problem}

\author{Sam Armon}
\address{University of Southern California, Los Angeles, CA 90089, USA}
\email{armon@usc.edu}

\author{Joshua P. Swanson}
\address{University of Southern California, Los Angeles, CA 90089, USA}
\email{swansonj@usc.edu}

\dedicatory{Dedicated to the memory of Adriano Garsia}
\date{\today}

\newtheorem{thm}{Theorem}
\numberwithin{thm}{section}
\newtheorem{lemma}[thm]{Lemma}
\newtheorem{prop}[thm]{Proposition}
\newtheorem{cor}[thm]{Corollary}

\theoremstyle{definition}
\newtheorem{definition}[thm]{Definition}
\newtheorem{example}[thm]{Example}
\newtheorem{prob}[thm]{Problem}
\newtheorem{remark}[thm]{Remark}

\DeclareMathOperator{\gr}{gr}
\DeclareMathOperator{\Hom}{Hom}

\newcommand\CC{\mathbb{C}}
\newcommand\QQ{\mathbb{Q}}
\newcommand\ZZ{\mathbb{Z}}

\newcommand{\fgtil}{\widetilde{\mathfrak{g}}}

\newcommand\ch{\Ch}
\newcommand\comaj{\mathrm{comaj}}
\newcommand\Des{\mathrm{Des}}
\newcommand\GL{\mathrm{GL}}
\newcommand\grch{\mathrm{GrCh}}
\newcommand\maj{\mathrm{maj}}
\newcommand\Neg{\mathrm{Neg}}
\newcommand\negg{\mathrm{neg}}

\newcommand\SYT{\mathrm{SYT}}
\newcommand\Srm{S}
\newcommand\Trm{\mathrm{T}}

\newcommand\xb{\mathbf{x}}
\newcommand\yb{\mathbf{y}}

\newcommand\Acal{\mathcal{A}}
\newcommand\Atil{\widetilde{A}}
\newcommand\Lcal{\mathcal{L}}
\newcommand\Ltil{\widetilde{\mathcal{L}}}
\newcommand\Qtil{\widetilde{Q}}
\newcommand\ptil{\tilde{p}}
\newcommand\stil{\tilde{s}}
\newcommand\Stil{\widetilde{S}}
\newcommand\Util{\widetilde{U}}
\newcommand\Trmtil{\widetilde{\Trm}}

\newcommand\Ch{\mathrm{Ch}}
\newcommand\FrobCh{\mathrm{FrobCh}}

\usepackage{xcolor}
\usepackage[colorinlistoftodos]{todonotes}

\begin{document}

\begin{abstract}
    Thrall's problem asks for the Schur decomposition of the higher Lie modules $\mathcal{L}_\lambda$, which are defined using the free Lie algebra and decompose the tensor algebra as a general linear group module. Although special cases have been solved, Thrall's problem remains open in general. We generalize Thrall's problem to the free Lie superalgebra, and prove extensions of three known results in this setting: Brandt's formula, Klyachko's identification of the Schur--Weyl dual of $\mathcal{L}_n$, and Kr{\'a}skiewicz--Weyman's formula for the Schur decomposition of $\mathcal{L}_n$. The latter involves a new version of the major index on super tableaux, which we show corresponds to a $q,t$-hook formula of Macdonald.
\end{abstract}

\maketitle

\section{Introduction}
\subsection{Super Thrall's problem}
Thrall \cite{Thrall} famously introduced a certain canonical decomposition of the tensor algebra coming from free Lie algebras:
\begin{equation}\label{eq:Thrall}
    \Trm(V) \cong \bigoplus_\lambda \Lcal_\lambda(V).
\end{equation}
Here $V = \CC^N$, $\lambda$ ranges over all integer partitions, the $\Lcal_\lambda = \Lcal_\lambda(V)$ are called \textit{Lie modules}, and the isomorphism is as $\GL(V)$-modules. See \Cref{sec:Thrall-classical} for details. \textit{Thrall's problem} is to explicitly decompose the $\Lcal_\lambda$ in terms of irreducible representations. Equivalently, we seek the Schur decomposition of the characters $\ch(\Lcal_\lambda; \xb)$ in the limit $N \to \infty$, where $\xb = x_1, x_2, \ldots$. Thrall's problem remains open in general, and has received significant attention in the literature since its formulation \cite{Brandt,Kly,BBG,Gar,GR,Sun,KW,Sch,AS}. See Vic Reiner's lecture slides \cite{Rei} for a lucid introduction to Thrall's problem and Adriano Garsia's lecture notes \cite{Gar} for a detailed introduction to the theory of free Lie algebras and many of the results we extend.

We study a supersymmetric extension of Thrall's problem. In \Cref{thm:super-Thrall}, we extend \eqref{eq:Thrall} to a canonical decomposition of the tensor superalgebra coming from free Lie superalgebras:
\begin{equation}\label{eq:super-Thrall}
    \Trmtil(V) \cong \bigoplus_A \widetilde{\Lcal}_A(V).
\end{equation}
Here $V = V_0 \oplus V_1 = \CC^N \oplus \CC^M$ is a $\ZZ/2$-graded vector space (a \textit{super vector space}), $A$ ranges over all infinite $\ZZ_{\geq 0}$-valued matrices $A=(a_{i,j})_{i,j \geq 0}$ with finite support and $a_{0,0}=0$, the $\Ltil_A = \Ltil_A(V)$ are called \textit{super Lie modules}, and the isomorphism is as $\GL(\CC^N) \oplus \GL(\CC^M)$-modules. See \Cref{sec:super-Thrall} for details. We consider the problem of explicitly decomposing $\Ltil_A$ into irreducible representations when $N=M \to \infty$ under the diagonal action $\GL(\CC^N) \hookrightarrow \GL(\CC^N) \oplus \GL(\CC^N)$. Equivalently, we seek the Schur decomposition of the characters $\ch(\Ltil_A; \xb, \xb)$. We recover the classical case as $\Ltil_A = \Lcal_\lambda$ when $a_{i,0}$ is the number of copies of $i$ in $\lambda$, and $a_{i,j} = 0$ for $j > 0$.

\subsection{A key special case}
The most important component of \eqref{eq:Thrall} is $\lambda = (n)$, where $\Lcal_n = \Lcal_{(n)}$ is the degree-$n$ homogeneous component of the free Lie algebra of $V$ in $\Trm(V)$. The decomposition of the higher Lie modules $\Lcal_\lambda$ may be determined by $\ch(\Lcal_n)$ from the Littlewood--Richardson rule and plethysm.

In the super case, $\widetilde{\Lcal}_{n,m}$ is the bidegree-$(n,m)$ homogeneous component of the free Lie superalgebra of $V_0 \oplus V_1$ in $\Trmtil(V_0 \oplus V_1)$. We have $\widetilde{\Lcal}_{n,0} = \Lcal_n$ and $\widetilde{\Lcal}_{n,m} = \widetilde{\Lcal}_A$ where $a_{n,m}=1$ and $a_{i,j} = 0$ otherwise. As before, the decomposition of the higher super Lie modules $\widetilde{\Lcal}_A$ may be determined by $\ch(\widetilde{\Lcal}_{n,m})$ from the Littlewood--Richardson rule and plethysm.

\textit{Brandt's formula} \cite{Brandt} gives the decomposition of $\Lcal_n$ in the power sum basis. We prove the following generalization:
\begin{thm}\label{thm:super-brandt}
    The $\GL(\CC^N)$-character of $\Ltil_{n,m}$ is given by
    \begin{equation}
        \ch(\Ltil_{n,m}; \xb_N) = \frac{1}{n+m} \sum_{d \mid \gcd(n,m)} (-1)^{m + \frac{m}{d}} \mu(d) \binom{\frac{n+m}{d}}{\frac{m}{d}} p_d(\xb_N)^\frac{n+m}{d}.
    \end{equation}
\end{thm}
\noindent The $m=0$ case of \Cref{thm:super-brandt} reduces to Brandt's formula. We also obtain the dimension of $\Ltil_{n, m}$, generalizing \textit{Witt's formula} \cite{Witt} in the $m=0$ case:

\begin{cor}[{\cite{Pet}}]\label{cor:super-Witt}
    When $V_0 = V_1 = \mathbb{C}^N$ and $n+m > 0$,
      \[ \dim \Ltil_{n,m} = \frac{1}{n+m} \sum_{d \mid \gcd(n, m)} (-1)^{m + \frac{m}{d}} \mu(d) \binom{\frac{n+m}{d}}{\frac{m}{d}} N^{\frac{n+m}{d}}. \]
\end{cor}

Klyachko \cite{Kly} described the Schur--Weyl dual of $\Lcal_n$ in terms of certain induced representations. To state our generalization, we need some notation. Let $C_{n+m}$ be the cyclic subgroup of the symmetric group $S_{n+m}$ generated by the long cycle $\pi_{n+m} = (1\,2\,\cdots\,n+m)$. Let $\chi^r \colon C_{n+m} \to \CC$ be the character with $\chi^r(\pi_{n+m}) = \exp(2\pi ir/(n+m))$. Finally, let $\chi^\mathrm{cyc}$ be the character of the action where $\pi_{n+m}$ cyclically increments $m$-element subsets of $[n+m]$. We prove the following: 
\begin{thm}\label{thm:super-induced}
We have
\begin{equation}
    \Ch(\Ltil_{n,m}) =
    \begin{cases}
        \FrobCh((\chi^\mathrm{cyc} \otimes \chi^1)\uparrow_{C_{n+m}}^{S_{n+m}}) & \text{if $m$ is odd} \\
        \FrobCh((\chi^\mathrm{cyc} \otimes \chi^{m/2+1})\uparrow_{C_{n+m}}^{S_{n+m}}) & \text{if $m$ is even.}
    \end{cases}
\end{equation}
\end{thm}
\noindent When $m=0$, the $\chi^\mathrm{cyc}$ factor may be omitted, recovering the $\Lcal_n$ case.

\subsection{Super Schur decompositions}

Kra{\'s}kiewicz--Weyman determined the Schur decomposition of $\Lcal_n$ \cite{KW}. We generalize this result to $\widetilde{\Lcal}_{n,m}$. Our argument introduces several super extensions of existing concepts in tableau combinatorics.

Let $\SYT_\pm(\lambda)$ denote the set of standard Young tableaux of shape $\lambda \vdash n$ where each cell is marked as \textit{positive} or \textit{negative}. For $\mathcal{T} \in \SYT_\pm(\lambda)$, let $\Neg(\mathcal{T}) \subseteq [n]$ be the set of all $i$ such that the box containing $i$ in $\mathcal{T}$ is marked as negative. Write $\negg(\mathcal{T}) \coloneqq |\Neg(\mathcal{T})|$.
\begin{definition}
    The \textit{(super) descent set} of $\mathcal{T} \in \SYT_\pm(\lambda)$ is the set $\Des(\mathcal{T})$ of all $i=1, \ldots, n-1$ such that either $i+1 \not\in \Neg(\mathcal{T})$ and $i+1$ appears in a strictly lower row of $\mathcal{T}$ than $i$, or $i \in \Neg(\mathcal{T})$ and $i+1$ does not appear in a strictly lower row of $\mathcal{T}$ than $i$. The \textit{(super) major index} of $\mathcal{T} \in \SYT_\pm(\lambda)$ is
\begin{equation}\label{eq:super-maj}
    \maj(\mathcal{T}) = \sum_{i \in \Des(\mathcal{T})} i.
\end{equation}
\end{definition}
\noindent Here we use English notation, with the longest row on top, and $\Neg(\mathcal{T}) = \varnothing$ recovers the usual descent set and major index. 

Our super generalization of Kra{\'s}kiewicz--Weyman's theorem is as follows. When $m=0$, we recover Kra{\'s}kiewicz--Weyman's original result \cite{KW} in type $A$.
\begin{thm}\label{thm:super-KW}
    The multiplicity of the Schur module $V^\lambda$ in $\Ltil_{n,m}$ is given by
    \begin{equation}\label{eq:super-KW}
        |\{ \mathcal{T} \in \SYT_\pm(\lambda) : \maj(\mathcal{T}) \equiv_{n+m} 1, \negg(\mathcal{T}) = m \}|.
    \end{equation}
\end{thm}

Our proof of \Cref{thm:super-KW} involves super analogues of the theory of $P$-partitions and symmetric function identities. For instance, we prove the following $q,t$-hook formula, which was announced in \cite{BS}:
\begin{thm}\label{thm:maj-neg-hook}
    For any $\lambda \vdash n$,
    \begin{equation}\label{eq:maj-neg-hook}
        \sum_{\mathcal{T} \in \SYT_\pm(\lambda)} q^{\maj(\mathcal{T})} t^{\negg(\mathcal{T})} = [n]_q! \prod_{(r, c) \in \lambda} \frac{q^{r-1} + tq^{c-1}}{[h(r, c)]_q}
    \end{equation}
    where $[n]_q \coloneqq 1 + q + \cdots + q^{n-1}$ and $[n]_q! \coloneqq [1]_q[2]_q \cdots [n]_q$.
\end{thm}

By contrast, Kra{\'s}kiewicz--Weyman's proof uses a character computation to identify $\Lcal_n$ as isomorphic to a submodule of the classical coinvariant algebra. See \cite[\S 8]{Gar} for a detailed discussion. The coinvariant algebra is very well-studied and intimately related to geometry and topology. It would be interesting to extend this connection, which we leave as an open problem.
\begin{prob}
    Find an interpretation of \Cref{thm:super-KW} which directly relates super Lie modules and coinvariant algebras.
\end{prob}

\subsection{Paper organization}
 The rest of the paper is organized as follows. In \Cref{sec:background} we recall some necessary background on (super)symmetric functions and representation theory, and introduce the classical case of Thrall's problem. In \Cref{sec:Thrall} we define (free) Lie superalgebras and prove the super analogue of Brandt's formula (\Cref{thm:super-brandt}), before explicitly defining Thrall's problem for free Lie superalgebras. In \Cref{sec:super-induced} we determine the Schur--Weyl dual of $\Ltil_{n,m}$, proving \Cref{thm:super-induced}. In \Cref{sec:super-tab} we develop some super tableau combinatorics and prove the $q,t$-hook formula (\Cref{thm:maj-neg-hook}), which is in turn used to prove \Cref{thm:super-KW} in \Cref{sec:super-KW}. Further directions are discussed in \Cref{sec:addl}.

\subsection{Acknowledgments}

Swanson was partially supported by NSF DMS-2348843. The authors would like to thank Sheila Sundaram for the question which led to the present paper.

\section{Background}\label{sec:background}
\subsection{Symmetric functions}\label{sec:symf}
We first review some background on tableau combinatorics and symmetric functions, and fix notation. For full details, see \cite[Chapter~7]{Stan}. An \emph{integer partition} $\lambda$ is a weakly decreasing sequence of positive integers $\lambda = (\lambda_1 \ge \cdots \ge \lambda_k > 0)$, and we let $\mathrm{Par}$ denote the set of all integer partitions. If $\sum_i \lambda_i = n$, then we denote this by $\lambda \vdash n$. By a slight abuse of notation we identify $\lambda$ with its \emph{Ferrers diagram}, which is the subset of $\ZZ_{>0} \times \ZZ_{>0}$ consisting of $\lambda_i$ left-justified cells in row $i$. We draw our partitions using English notation, so that the longest row appears at the top. The \emph{length} of $\lambda$, denoted $\ell(\lambda)$, is the number of nonempty rows in its Ferrers diagram. We sometimes will also use the notation $\lambda = (1^{a_1}2^{a_2} \cdots)$ to denote the partition with $a_i$ parts equal to $i$. 

Assume $\lambda \vdash n$, and let $[n] = \{ 1, 2, \ldots, n \}$. A \emph{standard tableau} of shape $\lambda$ is a bijective map $T : \lambda \to [n]$ that is strictly increasing along the rows and columns of $\lambda$. The set of all standard tableaux of shape $\lambda$ will be denoted by $\SYT(\lambda)$. For $i = 1, \ldots, n-1$, we say that $i$ is a \emph{descent} of $T$ if $i+1$ appears in a strictly lower row of $T$ than $i$. Let
\[
\Des(T) \coloneqq \{ i : i \text{ is a descent of } T \} \subseteq [n-1].
\]
\begin{example}
    The standard tableau
    \[
    T = 
    \ytableausetup{boxsize=1.5em,centertableaux}
    \begin{ytableau}
        1 & 3 & 4 & 6 \\
        2 & 5 \\
        7
    \end{ytableau}
    \in \SYT(4,2,1)
    \]
    has $\Des(T) = \{ 1, 4, 6 \}$.
\end{example}
Let $\xb = (x_1,x_2,\ldots)$. For $n \ge 2$ and $D \subseteq [n-1]$, the \emph{fundamental quasisymmetric function} $Q_{n,D}(\xb)$ is given by
\[
Q_{n,D}(\xb) = \sum_{\substack{a_1 \le a_2 \le \cdots \le a_n, \\ a_i = a_{i+1} \Rightarrow i \not\in D}} x_{a_1} x_{a_2} \cdots x_{a_n}.
\]

The \emph{Schur function} $s_\lambda(\xb)$ indexed by $\lambda$ is written in terms of quasisymmetric functions as follows:
\[
s_\lambda(\xb) = \sum_{T \in \SYT(\lambda)} Q_{n,\Des(T)}(\xb).
\]
Schur functions form a basis for the ring of \emph{symmetric functions} $\Lambda(\xb) \subseteq \ZZ[[x_1,x_2,\ldots]]$, the subspace consisting of all power series of bounded degree which are unchanged by permuting variable indices. Another important basis for the ring of symmetric functions is formed by the \emph{power sum} symmetric functions $p_\lambda(\xb)$, which are given by
\[
p_d(\xb) = x_1^d + x_2^d + \cdots, \qquad p_\lambda(\xb) = p_{\lambda_1}(\xb)p_{\lambda_2}(\xb) \cdots
\]
for any $d \ge 1$ and any $\lambda = (\lambda_1 \ge \lambda_2 \ge \cdots > 0)$.

\subsection{Representation theory of $\GL(\CC^N)$}
We next provide some background about the the representation theory of the general linear group $\GL(\CC^N)$ of invertible linear transformations of $\CC^N$. We discuss its Schur--Weyl duality with the symmetric group $S_k$ --- a complete exposition of which can be found, for instance, in \cite{Ful} --- and then recount the relevant facts about its character theory and (bi)graded representations.

Let $N \ge 1$. A \emph{$\GL(\CC^N)$-representation} $E$ is a $\CC$-vector space equipped with an action of $\GL(\CC^N)$. All representations in this paper will be polynomial. The irreducible representations of $\GL(\CC^N)$ are the \emph{Schur modules} $V^\lambda$, and are indexed by integer partitions $\lambda = (\lambda_1 \ge \cdots \ge \lambda_N \ge 0)$ with at most $N$ non-zero parts. 

By the Schur--Weyl duality between the representation theory of $\GL(\CC^N)$ and the symmetric group $S_k$, the irreducible representations $V^\lambda$ may be constructed as follows. The irreducible representations of $S_k$ are the \emph{Specht modules} $S^\lambda$, which are indexed by partitions $\lambda \vdash k$. The \emph{Frobenius characteristic map} $\FrobCh$ is an isomorphism between the space of class functions on $S_k$ and the space of homogeneous degree $k$ symmetric functions, defined by $\FrobCh(S^\lambda) = s_\lambda(\xb)$ for any $\lambda \vdash k$. Letting $V = \CC^N$, the space $V^{\otimes k}$ admits a natural $\GL(V) \times S_k$-action, where $\GL(V)$ acts diagonally on $V^{\otimes k}$ on the left, and $S_k$ acts on the right by permuting tensor factors. For any $S_k$-module $M$, the \emph{Schur--Weyl dual} of $M$ is the $\GL(V)$-module
\[
E(M) = V^{\otimes k} \otimes_{\CC S_k} M.
\]
Then the irreducible representations of $\GL(V)$ are given as the Schur--Weyl duals of all $S^\lambda$ for $\ell(\lambda) \le N$.

Upon fixing a basis for $\CC^N$, we may identify $\GL(\CC^N)$ with the group of invertible $N \times N$ matrices over $\CC$. The \emph{character} of a $\GL(\CC^N)$-representation $E$ is the trace of the action of the diagonal matrix $\mathrm{diag}(x_1,\ldots,x_N)$ on $E$, denoted $\Ch(E;\xb_N)$. The character $\Ch(E;\xb_N)$ is a polynomial in the variables $\xb_N = x_1,\ldots,x_N$, and the character of the Schur module $V^\lambda$ is
\[
\Ch(V^\lambda;\xb_N) = s_\lambda(x_1,\ldots,x_N,0,0,\ldots) \in \ZZ[x_1,\ldots,x_N],
\]
the Schur polynomial in $N$ variables. We then have that
\[
\lim_{N \to \infty} \Ch(V^\lambda;\xb_N) = s_\lambda(\xb) = \FrobCh(S^\lambda),
\]
so in practice we will often let $N \to \infty$ and omit reference to the underlying vector space. By the $\ZZ$-linear independence of Schur functions, computing the irreducible decomposition of a $\GL(V)$- or $S_k$-module is equivalent to computing the expansion of its character (resp. Frobenius characterstic) in the Schur basis.

A \emph{graded} (resp. \emph{bigraded}) $\GL(\CC^N)$-representation is a graded (resp. doubly graded) $\CC$-vector space $E = \bigoplus_n E_n$ (resp. $F = \bigoplus_{n,m} F_{n,m}$) equipped with a $\GL(\CC^N)$-action, such that the graded pieces $E_n$ (resp. $F_{n,m}$) are themselves $\GL(\CC^N)$-representations. The \emph{graded characters} $\grch(E;\xb;q)$ and $\grch(F;\xb,q,t)$ of $E$ and $F$ (resp.) are defined as
\[
\grch(E;\xb_N,q) = \sum_n \Ch(E_n;\xb_N)q^n, \quad \grch(F;\xb_N,q,t) = \sum_{n,m} \Ch(F_{n,m};\xb_N)q^nt^m.
\]

We will later consider representations of $\GL(\CC^N) \oplus \GL(\CC^M)$. After fixing bases for $\CC^N$ and $\CC^M$, we may consider the action of $\mathrm{diag}(x_1,\ldots,x_N,y_1,\ldots,y_M)$ on a $\GL(\CC^N) \oplus \GL(\CC^M)$-representation $E$. As above, the \emph{character} $\Ch(E;\xb_N,\yb_M)$ of $E$ is the trace of this action on $E$, which is a polynomial in $\CC[x_1,\ldots,x_N;y_1,\ldots,y_M]$, symmetric in the $x$- and $y$-variables separately. If $E = \bigoplus_{n,m} E_{n,m}$ is a bigraded $\GL(\CC^N) \oplus \GL(\CC^M)$-representation, then its graded character is defined as
\[
\grch(E;\xb_N,\yb_M,q,t) = \sum_{n,m} \Ch(E_{n,m};\xb_N,\yb_M)q^nt^m.
\]
\subsection{The classical case of Thrall's problem}\label{sec:Thrall-classical}
We now define Thrall's problem and recall some of the progress that has been made on it. The tensor algebra decomposition considered by Thrall is described in \eqref{PBW-decomp}, and the higher Lie modules $\Lcal_\lambda$ are defined in \Cref{def:lie-mod}. We then outline the three results related to Thrall's problem which we seek to generalize: Klaychko's identification of the Schur--Weyl dual of $\Lcal_n$ (\Cref{thm:kly}), Brandt's formula (\Cref{thm:brandt}), and Kra{\'s}kiewicz--Weyman's Schur decomposition of $\Ch(\Lcal_n;\xb)$ (\Cref{thm:KW}).

Let $V$ be a $\CC$-vector space with basis $X = \{ x_1,\ldots,x_N \}$, and let $\Trm(V) = \bigoplus_{n \ge 0} V^{\otimes n}$ denote its tensor algebra, which is naturally a graded $\GL(\CC^N)$-module. The \emph{free Lie algebra} $\Lcal(V)$ is the Lie subalgebra of $\Trm(V)$ generated by $X$. Then $\Lcal(V) = \bigoplus_{n \ge 1} \Lcal_n(V)$ inherits a grading from $\Trm(V)$, where $\Lcal_n(V) \coloneqq \Lcal(V) \cap V^{\otimes n}$, endowing $\Lcal(V)$ with the structure of a graded $\GL(\CC^N)$-module. See \cite{Reu} for more background on the free Lie algebra.

It is well-known that the universal enveloping algebra $U(\Lcal(V))$ of $\Lcal(V)$ is naturally isomorphic to the tensor algebra $\Trm(V)$. On the other hand, the Poincar\'e--Birkhoff--Witt theorem gives $U(\Lcal(V)) \cong S(\Lcal(V))$, where $\Srm(W) = \bigoplus_{m \ge 0} \Srm^m(W)$ denotes the symmetric algebra of a complex vector space $W$. Hence we have the following decomposition of $\Trm(V)$ as a $\GL(\CC^N)$-module:
\begin{align}
\Trm(V) &\cong U(\Lcal(V) \cong \Srm(\Lcal(V)) \cong \bigotimes_{n \ge 1} \Srm(\Lcal_n(V)) = \bigotimes_{n \ge 1} \left( \bigoplus_{m \ge 0} \Srm^m(\Lcal_n(V)) \right) \nonumber \\
&\cong \bigoplus_{\lambda = (1^{m_1}2^{m_2} \cdots)} \Srm^{m_1}(\Lcal_1(V)) \otimes \Srm^{m_2}(\Lcal_2(V)) \otimes \cdots. \label{PBW-decomp}
\end{align}
We more carefully derive an analogous decomposition of the tensor algebra $\Trmtil(V)$ of a super vector space in \Cref{thm:super-Thrall}. The isomorphism in \eqref{PBW-decomp} motivates the following definition.
\begin{definition}\label{def:lie-mod}
For $\lambda = (1^{m_1}2^{m_2} \cdots)$, the \emph{higher Lie module} $\Lcal_\lambda(V)$ is:
\[
\Lcal_\lambda(V) = \Srm^{m_1}(\Lcal_1(V)) \otimes \Srm^{m_2}(\Lcal_2(V)) \otimes \cdots.
\]
\end{definition}
Thus the higher Lie modules $\{ \Lcal_\lambda(V) : \lambda \in \mathrm{Par} \}$ yield a decomposition of $\Trm(V)$ as a $\GL(\CC^N)$-module, and it is therefore natural to ask for their irreducible decompositions. This problem was originally posed by Thrall \cite{Thrall}. By taking characters, an equivalent formulation of Thrall's problem is as follows.
\begin{prob}[Thrall's problem]
  Determine the coefficients $a_\mu \in \ZZ_{\ge 0}$ in the Schur expansion of $\Ch(\Lcal_\lambda(V);\xb_N)$:
  \[
  \Ch(\Lcal_\lambda(V);\xb_N) = \sum_\mu a_\mu s_\mu(\xb_N).
  \]
\end{prob}
\noindent As noted above, we may take $N \to \infty$ and work instead with symmetric functions; in this case we will simply write $\Lcal_\lambda$ for the higher Lie module, omitting the reference to the underlying vector space.

By \eqref{PBW-decomp} and the Littlewood--Richardson rule, it would suffice to determine the irreducible components of $\Lcal_{(a^b)} = S^b(\Lcal_a)$ for the purposes of Thrall's problem. On the level of characters, we have
\[
\Ch(\Lcal_{(a^b)};\xb_N) = h_b[\Ch(\Lcal_a)],
\]
where the right-hand side denotes the plethysm of $h_b$ with $\Ch(\Lcal_a)$ (see, for instance, \cite[Appendix~2]{Stan}). Thrall's problem remains open in general, although the single-row case follows from work of Klyachko \cite{Kly} and Kra{\'s}kiewicz--Weyman \cite{KW} by identifying the Schur--Weyl dual of $\Lcal_n$.

Let $\pi_n = (12 \cdots n) \in S_n$, and let $\omega_n$ denote a primitive $n$-th root of unity. Let $C_n = \langle \pi_n \rangle \le S_n$ denote the cyclic subgroup of order $n$ generated by $\pi_n$. The irreducible characters $\chi^1, \ldots, \chi^n$ of $C_n$ are given by $\chi^k(\pi_n) = \omega_n^k$. Klyachko proved that the Schur--Weyl dual of $\Lcal_n$ is obtained by inducing the representation $\chi^1$ of $C_n$ up to $S_n$:
\begin{thm}[\cite{Kly}]\label{thm:kly}
For any $n \ge 1$, the Schur--Weyl dual of $\Lcal_n$ is $\chi^1 \uparrow_{C_n}^{S_n}$; that is,
\[
\Ch(\Lcal_n;\xb) = \FrobCh(\chi^1 \uparrow_{C_n}^{S_n};\xb).
\]
\end{thm}
This result is equivalent to an expansion for $\Ch(\Lcal_n;\xb)$ in the power sum basis, which was also found directly by Brandt:
\begin{thm}[\cite{Brandt}]\label{thm:brandt}
For any $n \ge 1$,
\[
\Ch(\Lcal_n;\xb) = \frac{1}{n} \sum_{d \mid n} \mu(d) p_d(\xb)^{\frac{n}{d}},
\]
where $\mu(-)$ denotes the M{\"o}bius function.
\end{thm}
Kra{\'s}kiewicz--Weyman leverage the above expansion to determine the expansion of $\Ch(\Lcal_n;\xb)$ in the Schur basis.
\begin{thm}[\cite{KW}]\label{thm:KW}
For any $n \ge 1$, we have
\[
\Ch(\Lcal_n;\xb) = \sum_{\mu \vdash n} a_{\mu,1} s_\mu(\xb),
\]
where $a_{\mu,1} = |\{ T \in \SYT(\mu) : \maj(T) \equiv_n 1 \}|$.
\end{thm}
Outside of the single-row case, the only other known case of Thrall's problem is when $\lambda = (2^b)$, in which case $\Ch(\Lcal_{(2^b)})$ is given by a known plethystic identity:
\[
\Ch(\Lcal_{(2^b)}) = h_b[e_2] = \sum_{\mu} s_\mu,
\]
where the sum is over all $\mu \vdash 2b$ whose columns all have even length (see e.g. \cite[Ex.~I.8.6(b)]{Mac}).

\section{Super Thrall's problem}\label{sec:Thrall}
We now wish to extend Thrall's problem to the setting of free Lie superalgebras. In \Cref{sec:lie-super} we recall some facts about (free) Lie superalgebras, and then obtain a super generalization of Brandt's formula, \Cref{thm:super-brandt}. In \Cref{sec:super-Thrall} we identify the analogue of the higher Lie modules for the free Lie superalgebra, and state the super generalization of Thrall's problem in \Cref{thm:super-Thrall}.

\subsection{Free Lie superalgebras}\label{sec:lie-super}
Here we state the necessary background on Lie superalgebras, and then describe the structure of the free Lie superalgebra $\Ltil(V)$ as a bigraded $\GL(\CC^N) \oplus \GL(\CC^M)$-module. A result of Petrogradsky (\Cref{thm:free-super-lie-char}) determines the $\GL(\CC^N) \oplus \GL(\CC^M)$-character of $\Ltil(V)$, and we use this to obtain a character formula for the bigraded components of $\Ltil(V)$ in \Cref{thm:super-bi-brandt}, from which we obtain a proof of \Cref{thm:super-brandt}. We work over $\CC$ throughout.

A \textit{superalgebra} $\Atil$ is a $\ZZ/2$-graded vector space $\Atil = \Atil_0 \oplus \Atil_1$ with an associative, bilinear multiplication with a unit and satisfying $\Atil_i \Atil_j \subset \Atil_{i+j}$, where the indices are taken modulo $2$. For $x \in \Atil_i$, write $|x| = i$. The \textit{super commutator} on $\Atil$ is defined by extending
  \[ [x, y] \coloneqq xy - (-1)^{|x||y|} yx \]
bilinearly. It satisfies 
\begin{itemize}
\item[S1.] $[x,y] = -(-1)^{|x||y|}[y,x]$,
\item[S2.] $(-1)^{|x||z|}[x,[y,z]] + (-1)^{|z||y|}[z,[x,y]] + (-1)^{|y||x|}[y,[z,x]] = 0$.
\end{itemize}

\begin{example}
    The \textit{tensor superalgebra} of $V = V_0 \oplus V_1 = \CC^N \oplus \CC^M$ is
      \[ \Trmtil(V) \coloneqq \bigoplus_{d=0}^\infty (V_0 \oplus V_1)^{\otimes d} \]
    with the natural concatenation product. It has a bigrading where $\Trmtil_{n, m}(V)$ is spanned by tensors with $n$ factors from $V_0$ and $m$ factors from $V_1$. The $\ZZ/2$-grading is given by
      \[ \Trmtil_0(V) = \bigoplus_{n, m \geq 0} \Trmtil_{n, 2m}(V) \qquad\text{and}\qquad \Trmtil_1(V) = \bigoplus_{n, m \geq 0} \Trmtil_{n, 2m+1}(V). \]
    Both gradings respect the natural $\GL(\CC^N) \oplus \GL(\CC^M)$-action.
\end{example}

\begin{example}
    The \textit{symmetric superalgebra} of $V = V_0 \oplus V_1 = \CC^N \oplus \CC^M$ is
      \[ \Stil(V) \coloneqq S(V_0) \otimes \bigwedge(V_1) \]
    where $S(W)$ is the symmetric algebra and $\bigwedge(W)$ is the exterior algebra of the vector space $W$. As with $\Trmtil(V)$, it is a bigraded and $\ZZ/2$-graded $\GL(\CC^N) \oplus \GL(\CC^M)$-module. The super commutator is identically zero.
\end{example}

A \emph{Lie superalgebra} $\fgtil$ is a $\ZZ/2$-graded vector space $\fgtil = \fgtil_0 \oplus \fgtil_1$ with a bilinear operation $[-, -]$ satisfying (S1) and (S2) with $[\fgtil_i, \fgtil_j] \subset \fgtil_{i+j}$. As with all superalgebras, the tensor superalgebra $\Trmtil(V)$ and the symmetric superalgebra $\Stil(V)$ are Lie superalgebras under the super commutator.

\begin{example}
    The \textit{free Lie superalgebra} generated by $V = V_0 \oplus V_1 = \CC^N \oplus \CC^M$ is the Lie superalgebra $\Ltil(V)$ in $\Trmtil(V)$ generated by $V$. As with $\Trmtil(V)$, it is a bigraded and $\ZZ/2$-graded $\GL(\CC^N) \oplus \GL(\CC^M)$-module. Concretely, $\Ltil(V) \subset \Trmtil(V)$ is given by
    \begin{align*}
        \Ltil_{n, m}(V) &\coloneqq \text{span of iterated bracketings with $n$ terms from $V_0$, $m$ terms from $V_1$} \\
        \Ltil_0(V) &\coloneqq \text{span of iterated bracketings with evenly many terms from $V_1$} \\
        \Ltil_1(V) &\coloneqq \text{span of iterated bracketings with oddly many terms from $V_1$}.
    \end{align*}
    Note that $\Ltil_{0, 0}(V) = 0$, $\Ltil_{1, 0}(V) = V_0$, and $\Ltil_{0, 1}(V) = V_1$.
\end{example}

Petrogradsky \cite{Pet} found the doubly graded Hilbert series of $\Ltil(V)$, generalizing what is known in the classical case as \emph{Witt's formula}, which in turn yields a formula for the graded character of $\Ltil(V)$ as a $\GL(\CC^N) \oplus \GL(\CC^M)$-module.
\begin{thm}[\cite{Pet}]\label{thm:free-super-lie-char}
    The bihomogeneous $\GL(\CC^N) \oplus \GL(\CC^M)$-character of $\Ltil(V)$ is
    \[
    \grch(\Ltil(V);\xb_N,\yb_M,q,t) = -\sum_{d=1}^\infty \frac{\mu(d)}{d} \ln(1-(q^dp_d(\xb_N) - t^dp_d(-\yb_M))).
    \]
\end{thm}
The above result of Petrogradsky allows us to determine a character formula for the doubly graded pieces $\Ltil_{n,m}(V)$ of $\Ltil(V)$.
\begin{thm}\label{thm:super-bi-brandt}
    For $(n,m) \neq (0,0)$, the $\GL(\CC^N) \oplus \GL(\CC^M)$-character of $\Ltil_{n,m}(V)$ is
    \[
    \ch(\Ltil_{n,m}(V);\xb_N,\yb_M) = \frac{1}{n+m} \sum_{d \mid \gcd(m,n)} (-1)^\frac{m}{d} \mu(d) \binom{\frac{n+m}{d}}{\frac{m}{d}} p_d(\xb_N)^\frac{n}{d}p_d(-\yb_M)^\frac{m}{d}.
    \]
\end{thm}
\begin{proof}
    We have
    \begin{align*}
    \grch(\Ltil(V);\xb_N,\yb_M,q,t) &= -\sum_{d=1}^\infty \frac{\mu(d)}{d} \ln(1-(q^dp_d(\xb_N) - t^dp_d(-\yb_M))) \\
    &= \sum_{d=1}^\infty \frac{\mu(d)}{d} \sum_{s=1}^\infty \frac{(q^dp_d(\xb_N) - t^dp_d(-\yb_M))^s}{s} \\
    &= \sum_{d=1}^\infty \frac{\mu(d)}{d} \sum_{s=1}^\infty \frac{1}{s} \sum_{k=0}^s {\binom{s}{k}}(q^dp_d(\xb_N))^k(-t^dp_d(-\yb_M))^{s-k} \\
    &= \sum_{d=1}^\infty \sum_{\substack{0 \le k \le s < \infty, \\ (k,s) \neq (0,0)}} \frac{\mu(d)}{d} \frac{1}{s} {\binom{s}{k}} (q^dp_d(\xb_N))^k(-t^dp_d(-\yb_M))^{s-k} \\
    &= \sum_{d=1}^\infty \sum_{\substack{n : d \mid n, \\ m : d \mid m, \\ (n,m) \neq (0,0)}} \frac{\mu(d)}{n+m} \binom{\frac{n+m}{d}}{\frac{m}{d}} p_d(\xb_N)^\frac{n}{d}(-p_d(-\yb_M))^\frac{m}{d} q^nt^m.
    \end{align*}
    As $\ch(\Ltil_{n,m}(V);\xb_N,\yb_M) = [q^nt^m]\grch(\Ltil(V);\xb_N,\yb_M,q,t)$, the result follows.
\end{proof}
As a corollary to \Cref{thm:super-bi-brandt}, we obtain the super generalization of Brandt's formula, \Cref{thm:super-brandt}. Recall that we set $N=M$ and use the diagonal action here.
\begin{proof}[Proof (of \Cref{thm:super-brandt}).]
By \Cref{thm:super-bi-brandt}, we have
\begin{align*}
\ch(\Ltil_{n,m}(V);\xb_N) &= \frac{1}{n+m} \sum_{d \mid \gcd(m,n)} (-1)^\frac{m}{d} \mu(d) \binom{\frac{n+m}{d}}{\frac{m}{d}} p_d(\xb_N)^\frac{n}{d}p_d(-\xb_N)^\frac{m}{d} \\
&= \frac{1}{n+m} \sum_{d \mid \gcd(n,m)} (-1)^{m + \frac{m}{d}} \mu(d) \binom{\frac{n+m}{d}}{\frac{m}{d}} p_d(\xb_N)^\frac{n+m}{d}.
\end{align*}
\end{proof}

\subsection{Thrall's problem for free Lie superalgebras}\label{sec:super-Thrall}
We now define the super higher Lie modules (see \Cref{def:super-Lie}) which form a $\GL(\CC^N) \oplus \GL(\CC^M)$ decomposition of the tensor algebra of a super vector space described in \Cref{thm:super-Thrall}. We then state the super generalization of Thrall's problem.

The \textit{universal enveloping superalgebra} of a Lie superalgebra $\fgtil$ is intuitively the smallest superalgebra $\Util(\fgtil)$ which $\fgtil$ embeds into as a Lie superalgebra. Formally, $\Util(\fgtil)$ is a superalgebra with a map of Lie superalgebras $\iota \colon \fgtil \to \Util(\fgtil)$ satisfying the following universal property. For any superalgebra $\Atil$, composition with $\iota$ induces a bijection
  \[ \Hom_{\textrm{SuperLie}}(\fgtil, \Atil) \cong \Hom_{\textrm{SuperAlg}}(\Util(\fgtil), \Atil). \]
This leads to a natural construction of $\Util(\fgtil)$ as a quotient of $\Trmtil(\fgtil)$, which results in a filtration $\Util(\fgtil) = \bigcup_{d \geq 0} \Util_{\leq d}(\fgtil)$. The \textit{associated graded} superalgebra is $\gr \Util(\fgtil) \coloneqq \bigoplus_{d \geq 0} \Util_{\leq d}(\fgtil)/\Util_{\leq d-1}(\fgtil)$ where $\Util_{\leq -1}(\fgtil) \coloneqq 0$. Note that $\Util(\fgtil) \cong \gr \Util(\fgtil)$ as vector spaces.

\begin{example}\label{ex:U-of-free}
    The universal enveloping superalgebra of the free Lie superalgebra $\Ltil(V)$ is $\Trmtil(V)$. Indeed, if $\Atil$ is a superalgebra, $X$ is a basis for $V_0$, and $Y$ is a basis for $V_1$, we have natural identifications
    \begin{align*}
        \Hom_{\text{SuperAlg}}(\Util(\Ltil(V)), \Atil)
          &\cong \Hom_{\text{SuperLie}}(\Ltil(V), \Atil) \\
          &\cong \Hom_{\text{SuperSet}}(X \sqcup Y, \Atil_0 \sqcup \Atil_1) \\
          &\cong \Hom_{\text{SuperVec}}(V_0 \oplus V_1, \Atil_0 \oplus \Atil_1) \\
          &\cong \Hom_{\text{SuperAlg}}(\Trmtil(V), \Atil).
    \end{align*}
    Hence $\Hom_{\text{SuperAlg}}(\Util(\Ltil(V)), -) \cong \Hom_{\text{SuperAlg}}(\Trmtil(V), -)$, so by Yoneda's lemma,
      \[ \Util(\Ltil(V)) \cong \Trmtil(V). \]
\end{example}

The Poincar\'e--Birkhoff--Witt theorem for Lie superalgebras reads as follows. A standard consequence is that the map $\iota$ is an embedding, though we will only apply the result when $\fgtil = \Ltil(V)$, in which case injectivity is obvious by \Cref{ex:U-of-free}.

\begin{thm}[Super PBW]
   Let $\fgtil$ be a Lie superalgebra. Then there is a natural isomorphism of superalgebras
     \[ \Stil(\fgtil) \cong \gr \Util(\fgtil). \]
\end{thm}

We are nearly ready to state the super analogue of Thrall's decomposition of the tensor algebra. We first define analogues of Lie modules, which will arise naturally in the upcoming proof.

\begin{definition}\label{def:super-Lie}
    For a vector space $W$ and $j \in \ZZ_{\geq 0}$, let
    \[
    \Gamma_j(W) =
    \begin{cases}
        S(W) & \text{if $j$ is even} \\
        \bigwedge(W) & \text{if $j$ is odd}.
    \end{cases}
    \]
    Given a $\ZZ_{\ge 0}$-valued matrix $A = (a_{i,j})_{i,j \ge 0}$ with finite support and $a_{0,0} = 0$, the \emph{super Lie module} $\Ltil_A$ is:
    \begin{equation}\label{eq:super-higher-lie}
    \Ltil_A = \bigotimes_{i,j \ge 0} \Gamma_j^{a_{i,j}}(\Ltil_{i,j}).
    \end{equation}
    Note that for $j$ odd, the exterior power $\Gamma_j^a(W)$ is zero unless $a \leq \dim W$.
\end{definition}

\begin{thm}[Super Thrall decomposition]\label{thm:super-Thrall}
    Let $V = V_0 \oplus V_1 = \CC^N \oplus \CC^M$. Then
      \[ \Trmtil(V) = \bigoplus_A \Ltil_A \]
    as $\GL(\CC^N) \oplus \GL(\CC^M)$-modules, where the sum is over $\ZZ_{\geq 0}$-valued matrices $A = (a_{i,j})_{i,j \geq 0}$ with finite support and $a_{0,0} = 0$.
\end{thm}

\begin{proof}
    By \Cref{ex:U-of-free}, $\Trmtil(V) = \Util(\Ltil(V))$. Hence $\gr \Util(\Ltil(V)) = \Util(\Ltil(V))$. The super PBW theorem and routine properties now gives
    \begin{align*}
        \Trmtil(V)
          &= \Util(\Ltil(V)) = \gr \Util(\Ltil(V)) = \Stil(\Ltil(V)) = S(\Ltil(V)_0) \otimes \bigwedge(\Ltil(V)_1) \\
          &= S\left(\bigoplus_{n,m} \Ltil_{n, 2m}\right) \otimes \bigwedge\left(\bigoplus_{n,m} \Ltil_{n, 2m+1}\right)
           = \bigotimes_{n,m} S\left(\Ltil_{n, 2m}\right) \otimes \bigotimes_{n, m} \bigwedge\left(\Ltil_{n, 2m+1}\right) \\
          &= \bigotimes_{i, j} \Gamma_j\left(\Ltil_{i, j}\right)
           = \bigotimes_{i, j} \left( \bigoplus_{a_{i, j} \geq 0} \Gamma_j^{a_{i, j}}\left(\Ltil_{i, j}\right) \right)
           = \bigoplus_{A = (a_{i, j} \geq 0)} \bigotimes_{i, j} \Gamma_j^{a_{i, j}}\left(\Ltil_{i, j}\right) \\
          &= \bigoplus_{A = (a_{i, j} \geq 0)} \Ltil_A.\qedhere
    \end{align*}
\end{proof}

Now, if $N = M$, then $\Ltil$ inherits the structure of a $\GL(\CC^N)$-module under the diagonal inclusion $\GL(\CC^N) \hookrightarrow \GL(\CC^N) \times \GL(\CC^N)$. Then Thrall's problem may be generalized to the free Lie superalgebra as follows.
\begin{prob}[Super Thrall's problem]
  For $A = (a_{i, j} \geq 0)_{i, j}$ and $\lambda \in \mathrm{Par}$, what is the multiplicity of the Schur module $V^\lambda$ in the super higher Lie module $\Ltil_A$ where $N=M \to \infty$?
\end{prob}
As in the classical case, it would be sufficient to determine these multiplicities for the modules $S^a(\Ltil_{n,2m})$ and $\bigwedge^b(\Ltil_{n,2m+1})$ for general $a,b,n,m \in \ZZ_{\ge 0}$, by the Littlewood--Richardson rule. As a result, the corresponding problem for the modules $\Ltil_{n,m}$ is of particular importance.

\section{Schur--Weyl duals as induced characters}\label{sec:super-induced}
In this section we prove \Cref{thm:super-induced}, which identifies the Schur--Weyl dual of $\Ltil_{n,m}$. We first describe a certain representation of $C_{n+m}$ and then prove that, upon inducing up to $S_{n+m}$, the Frobenius character of the resulting representation agrees with $\Ch(\Ltil_{n,m})$.

Recall from above that $C_r \le S_r$ is the cyclic group of order $r$, generated by the cycle $\pi_r = (1 2 \cdots r) \in S_r$, and $\omega_r$ denotes a primitive $r$-th root of unity. The irreducible characters of $C_r$ are $\chi^1, \ldots, \chi^r$, given by $\chi^d(\pi_r) = \omega_r^d$ for $1 \le d \le r$.

The group $C_r$ acts naturally on $[r]$ by cyclic rotation. Fixing $n,m \ge 0$ with $(n,m) \neq (0,0)$, let $\binom{[n+m]}{m} = \{ S \subseteq [n+m] : |S| = m \}$. Then the action of $C_{n+m}$ on $[n+m]$ extends naturally to an action on $\binom{[n+m]}{m}$:
\[
\pi_{n+m}^k \cdot \{ i_1, \ldots, i_m \} = \{ i_1 + k\pmod{n+m}, \ldots, i_m + k\pmod{n+m} \}
\]
for any $k \ge 1$. Let $\chi^\mathrm{cyc} : C_{n+m} \to \CC$ denote the character of this latter action, so that
\[
\chi^\mathrm{cyc}(\pi_{n+m}^k) = |\{ S \in \binom{[n+m]}{m} : \pi_{n+m}^k \cdot S = S \}|.
\]

\begin{lemma}\label{cyc-char}
For any $n,m \ge 0$ with $(n,m) \neq 0$ and $1 \le k \le n+m$, let $d = \frac{n+m}{\gcd(k,n+m)}$. Then
\[
\chi^\mathrm{cyc}(\pi_{n+m}^k) = \begin{cases}
\binom{\frac{n+m}{d}}{\frac{m}{d}} & \text{if } d \mid m \\
0 & \text{if } d \nmid m.
\end{cases}
\]
\end{lemma}
\begin{proof}
Note that $\pi_{n+m}^k$ is a product of $\gcd(k,n+m) = \frac{n+m}{d}$ disjoint cycles, each of length $d$. Let $[n+m] = C_1 \sqcup \cdots \sqcup C_{\frac{n+m}{d}}$ denote the partition of $[n+m]$ given by the decomposition of $\pi_{n+m}^k$ into disjoint cycles. Now, consider the action of $\pi_{n+m}^k$ on $\binom{[n+m]}{m}$. A set $S \subseteq [n+m]$ is fixed by $\pi_{n+m}^k$ if and only if, for each $1 \le j \le \gcd(k,n+m)$, either $S \cap C_j = C_j$ or $S \cap C_j = \varnothing$. Thus the $m$-subsets of $[n+m]$ fixed by $\pi_{n+m}^k$ may be constructed by choosing $\frac{m}{d}$ of the sets $C_1,\ldots,C_{\frac{n+m}{d}}$. This is impossible if $d \nmid m$, and the result follows.
\end{proof}

For any subgroup $H \le S_r$ and $\chi$ a representation of $H$, the Frobenius characteristic of $\chi \uparrow_H^{S_r}$ may be easily expressed in the power sum basis (see e.g. \cite[Theorem~13]{Swa}). Here the cycle type of $\pi_{n+m}^k$ depends only on $\gcd(k,n+m)$, giving:
\begin{lemma}\label{induced-rep-char}
For any representation $\chi$ of $C_r$, the Frobenius characteristic of $\chi \uparrow_{C_r}^{S_r}$ is given by:
\[
\FrobCh(\chi \uparrow_{C_r}^{S_r}) = \frac{1}{r} \sum_{d \mid r} \sum_{k : \gcd(k,r) = \frac{r}{d}} \chi(\pi_r^k) p_d(\mathbf{x})^{\frac{r}{d}}.
\]
\end{lemma}
We now proceed to the proof of \Cref{thm:super-induced}.
\begin{proof}[Proof (of \Cref{thm:super-induced}).]
Throughout, let $\pi = \pi_{n+m}$ and $\omega = \omega_{n+m}$. 

First assume $m$ is odd, so that for any $d \mid m$, $m+\frac{m}{d}$ is always even. Then by \Cref{induced-rep-char}, we have:
\[
\FrobCh((\chi^\mathrm{cyc} \otimes \chi^1)\uparrow_{C_{n+m}}^{S_{n+m}}) = \frac{1}{n+m} \sum_{d \mid n+m} \sum_{k : \gcd(k,n+m) = \frac{n+m}{d}} \chi^\mathrm{cyc}(\pi^k) \chi^1(\pi^k) p_d(\mathbf{x})^{\frac{n+m}{d}},
\]
which by \Cref{cyc-char} simplifies to:
\[
\FrobCh((\chi^\mathrm{cyc} \otimes \chi^1)\uparrow_{C_{n+m}}^{S_{n+m}}) = \frac{1}{n+m} \sum_{\substack{d \mid n+m, \\ d \mid m}} \binom{\frac{n+m}{d}}{\frac{m}{d}} \sum_{k : \gcd(k,n+m) = \frac{n+m}{d}} \chi^1(\pi^k) p_d(\mathbf{x})^{\frac{n+m}{d}}
\]
Note that $d \mid n+m$ and $d \mid m$ if and only if $d \mid \gcd(n,m)$. Furthermore, if $\gcd(k,n+m) = \frac{n+m}{d}$, then $\chi^1(\pi^k) = \omega^k$ is a primitive $d$-th root of unity, and all primitive $d$-th roots of unity are obtained in this way. The M{\"o}bius function $\mu(d)$ is given as the sum of the primitive $d$-th roots of unity. Thus
\[
\FrobCh((\chi^\mathrm{cyc} \otimes \chi^1)\uparrow_{C_{n+m}}^{S_{n+m}}) = \frac{1}{n+m} \sum_{\substack{d \mid \gcd(n,m)}} \binom{\frac{n+m}{d}}{\frac{m}{d}} \mu(d) p_d(\mathbf{x})^{\frac{n+m}{d}} = \Ch(\Ltil_{n,m}),
\]
where the last equality follows from \Cref{thm:super-brandt}.

Now assume $m$ is even. Again applying \Cref{cyc-char} and \Cref{induced-rep-char}, we have
\begin{align*}
\FrobCh((\chi^\mathrm{cyc} &\otimes \chi^{\frac{m}{2}+1})\uparrow_{C_{n+m}}^{S_{n+m}}) \\
&= \frac{1}{n+m} \sum_{d \mid \gcd(n,m)} \binom{\frac{n+m}{d}}{\frac{m}{d}} \sum_{k : \gcd(k,n+m) = \frac{n+m}{d}} \chi^{\frac{m}{2}+1}(\pi^k) p_d(\mathbf{x})^{\frac{n+m}{d}}.
\end{align*}
Observe for any $d \mid m$ that $\frac{m}{d}$ is even if and only if $d \mid \frac{m}{2}$, and recall from above that if $\gcd(k,n+m) = \frac{n+m}{d}$, then $\omega^k$ is a primitive $d$-th root of unity. Thus in the case that $\frac{m}{d}$ is even, we have $(\omega^k)^{\frac{m}{2}} = 1$ since $d \mid \frac{m}{2}$, which equals $(-1)^{m+\frac{m}{d}}$ since $m$ and $\frac{m}{d}$ are even. On the other hand, if $\frac{m}{d}$ is odd, then $d$ must necessarily be even since $m$ is, so that
\[
(\omega^k)^{\frac{m}{2}} = ((\omega^k)^{\frac{d}{2}})^{\frac{m}{d}} = (-1)^{\frac{m}{d}} = (-1)^{m + \frac{m}{d}} = -1.
\]
Hence
\[
\chi^{\frac{m}{2}+1}(\pi^k) = (\omega^k)^{\frac{m}{2}}\omega^k = (-1)^{m+\frac{m}{d}}\omega^k.
\]
It then follows that
\begin{align*}
\FrobCh((\chi^\mathrm{cyc} &\otimes \chi^{\frac{m}{2}+1})\uparrow_{C_{n+m}}^{S_{n+m}}) \\
&= \frac{1}{n+m} \sum_{d \mid \gcd(n,m)} (-1)^{m + \frac{m}{d}} \mu(d) \binom{\frac{n+m}{d}}{\frac{m}{d}} p_d(\xb)^\frac{n+m}{d} = \Ch(\Ltil_{n,m}).
\end{align*}
\end{proof}

\section{Super tableau combinatorics}\label{sec:super-tab}
In this section we prove \Cref{thm:maj-neg-hook}, which is the key ingredient in our proof of \Cref{thm:super-KW}. We first recall some preliminaries on supersymmetric functions in \Cref{sec:supersym-fxns}. In \Cref{sec:super-des}, we define a new major index statistic on super tableaux and then prove \Cref{prop:q-ps} and \Cref{thm:s-ps}, which exhibit how the principal specializations of the super quasisymmetric and super Schur functions admit elegant formulae in terms of this new statistic. \Cref{thm:maj-neg-hook} then follows as an easy consequence.
\subsection{Supersymmetric functions}\label{sec:supersym-fxns}
We begin by recalling some facts about supersymmetric and super quasisymmetric functions, following the presentation in \cite{HHLRU}. Let $\yb = (y_1,y_2,\ldots)$. A power series $f(\xb;\yb) \in \Lambda(\xb) \otimes_\ZZ \Lambda(\yb)$ is \emph{supersymmetric} if performing the substitution $x_1 = t, y_1 = -t$ into $f$ results in an expression which is independent of $t$. Schur functions and power sums both admit supersymmetric analogs. The \emph{super power sum symmetric function} $\ptil_\lambda(\xb;\yb)$ indexed by $\lambda$ is:
\[
\ptil_\lambda(\xb;\yb) = \prod_{j=1}^{\ell(\lambda)} (p_{\lambda_j}(\xb) + (-1)^{\lambda_j+1}p_{\lambda_j}(\yb)).
\]

The super Schur functions were originally introduced by Berele--Regev \cite{BReg} in their study of the general linear Lie superalgebra, and they may be defined as follows. Let $\Acal_+ = \{ 1, 2, \ldots \}$ and $\Acal_- = \{ \overline{1}, \overline{2}, \ldots \}$. Endow the alphabet $\Acal = \Acal_+ \sqcup \Acal_-$ with the following total order:
\[
\Acal = \{ 1 < \overline{1} < 2 < \overline{2} < \cdots \}.
\]
We define a projection $\Acal \to \Acal_+$ by $i \mapsto i$ and $\overline{i} \mapsto i$ for all $i \ge 1$. Note that if $a_1 \le \cdots \le a_n$ is a weakly increasing sequence in $\Acal$ and $a_j \mapsto b_j$ under this projection, then $b_j \le a_j$ for all $j$ and $b_1 \le \cdots \le b_n$.
 
\begin{definition}[\cite{HHLRU}, Eq.~(23)]
For $n \ge 2$ and $D \subseteq [n-1]$, the \emph{super quasisymmetric function} $\Qtil_{n,D}(\xb;\yb)$ is given by
\[
\Qtil_{n,D}(\xb;\yb) = \sum_{\substack{a_1 \le a_2 \le \cdots \le a_n, \\ a_i = a_{i+1} \in \Acal_+ \Rightarrow i \not\in D, \\ a_i = a_{i+1} \in \Acal_- \Rightarrow i \in D}} z_{a_1} z_{a_2} \cdots z_{a_n},
\]
where $a_1 \le \cdots \le a_n$ is a weakly increasing sequence in $\Acal$, and $z_a = x_a$ for $a \in \Acal_+, z_b = y_b$ for $b \in \Acal_-$.
\end{definition}
As in the classical case, the super Schur function $\stil_\lambda(\xb;\yb)$ is given as a sum of super quasisymmetric functions:
\begin{definition}[\cite{HHLRU}, Prop.~2.4.2]\label{prop:super-s-to-q}
    The \emph{super Schur function} $\stil_\lambda(\xb;\yb)$ is given in terms of super quasisymmetric functions by
    \[
    \stil_\lambda(\xb;\yb) = \sum_{T \in \SYT(\lambda)} \Qtil_{|\lambda|,\Des(T)}(\xb;\yb).
    \]
\end{definition}

Schur functions and power sum symmetric functions are related via the \emph{Cauchy identity}, which extends to the supersymmetric setting as follows.
\begin{thm}[\cite{BRem}, Cor.~10(a)]\label{lem:s-p}
    For all $n \ge 1$,
    \[
    \sum_{\lambda \vdash n} s_\lambda(\xb) \stil_\lambda(\xb;\yb) = \sum_{\lambda \vdash n} \frac{1}{z_\lambda} p_\lambda(\xb) \ptil_\lambda(\xb;\yb),
    \]
    where if $\lambda = (1^{a_1}2^{a_2} \cdots)$, then $z_\lambda \coloneqq 1^{a_1}2^{a_2} \cdots a_1!a_2! \cdots$.
\end{thm}

\subsection{A major index statistic on super tableaux}\label{sec:super-des} We now define a new major index statistic in \Cref{def:super-maj} on objects called \emph{standard super tableaux}, and show how particular specializations of super quasisymmetric and super Schur functions may be written in terms of this statistic in \Cref{prop:q-ps} and \Cref{thm:s-ps}.
\begin{definition}
    A \emph{standard super tableau} of shape $\lambda \vdash n$ is a map $\mathcal{T} : \lambda \to \Acal$ that is strictly increasing along the rows and columns of $\lambda$, and contains exactly one of $i$ or $\overline{i}$ for each $i = 1, 2, \ldots, n$. Let $\SYT_\pm(\lambda)$ denote the set of standard super tableaux of shape $\lambda$, so that $|\SYT_\pm(\lambda)| = 2^n|\SYT(\lambda)|$. For $\mathcal{T} \in \SYT_\pm(\lambda)$, we let $\Neg(\mathcal{T}) \coloneqq \{ i \in \Acal_+ : \overline{i} \in \mathcal{T} \}$ and $\negg(\mathcal{T}) \coloneqq |\Neg(\mathcal{T})|$.
\end{definition}
\begin{example}
The standard super tableau
    \[
    \mathcal{T} = 
    \ytableausetup{boxsize=1.5em,centertableaux}
    \begin{ytableau}
        1 & {\color{red}\overline{3}} & 4 & 6 \\
        {\color{red}\overline{2}} & 5 \\
        {\color{red}\overline{7}}
    \end{ytableau}
    \in \SYT_\pm(4,2,1)
    \]
    has $\Neg(\mathcal{T}) = \{ 2, 3, 7 \}$ and $\negg(\mathcal{T}) = 3$.
\end{example}
Using this notion, we now define an appropriate generalization of the descent set for $\SYT_\pm(\lambda)$.
\begin{definition}
For $\lambda \vdash n$ and $\mathcal{T} \in \SYT_\pm(\lambda)$, let $\mathcal{T}_+ \in \SYT(\lambda)$ denote the image of $\mathcal{T}$ under the projection $\Acal \to \Acal_+$. For $i = 1, \ldots, n-1$, we say that $i$ is a \emph{super descent} of $\mathcal{T}$ if either 
\[
i \in \Des(\mathcal{T}_+) \text{ and } i+1 \not\in \Neg(\mathcal{T}), \text{ or } i \not\in \Des(\mathcal{T}_+) \text{ and } i \in \Neg(\mathcal{T}).
\]
Define
\[
\Des(\mathcal{T}) = \{ i : i \text{ is a super descent of } \mathcal{T} \} \subseteq [n-1].
\]
\end{definition}
Note that if $\Neg(\mathcal{T}) = \varnothing$, then the super descents of $\mathcal{T}$ are precisely the usual descents defined in Section \ref{sec:symf}.

\begin{definition}\label{def:super-maj}
For $D \subseteq [n-1]$ and $S \subseteq [n]$, define the \emph{relative major index} and the \emph{relative comajor index} respectively by
\[
\maj(D,S) \coloneqq \sum_{\substack{1 \le i \le n-1, \\ i \in D, i+1 \not\in S \\ \text{or } i \not\in D, i \in S}} i, \qquad \comaj(D,S) \coloneqq \sum_{\substack{1 \le i \le n-1, \\ i \in D, i+1 \not\in S \\ \text{or } i \not\in D, i \in S}} (n-i).
\]
For $\mathcal{T} \in \SYT_\pm(\lambda)$, we define the relative (co)major index by
\[
\maj(\mathcal{T}) \coloneqq \sum_{i \in \Des(\mathcal{T})} i, \quad \comaj(\mathcal{T}) \coloneqq \sum_{i \in \Des(\mathcal{T})} (n-i).
\]
\end{definition}
It follows from the above definition that
\[
\maj(\mathcal{T}) = \maj(\Des(\mathcal{T}_+),\Neg(\mathcal{T})), \quad \comaj(\mathcal{T}) = \comaj(\Des(\mathcal{T}_+),\Neg(\mathcal{T})).
\]
Note also that setting $S = \varnothing$ recovers the classical notions of major and comajor index on standard tableaux.
\begin{example}
If
    \[
    \mathcal{T} = 
    \ytableausetup{boxsize=1.5em,centertableaux}
    \begin{ytableau}
        1 & {\color{red}\overline{3}} & 4 & 6 \\
        {\color{red}\overline{2}} & 5 \\
        {\color{red}\overline{7}}
    \end{ytableau}
    \in \SYT_\pm(4,2,1) \quad \text{then} \quad
    \mathcal{T}_+ = 
    \ytableausetup{boxsize=1.5em,centertableaux}
    \begin{ytableau}
        1 & 3 & 4 & 6 \\
        2 & 5 \\
        7
    \end{ytableau}
    \in \SYT(4,2,1),
    \]
    so $\Des(\mathcal{T}_+) = \{ 1, 4, 6 \}$ and $\Neg(\mathcal{T}) = \{ 2, 3, 7 \}$. Therefore $\Des(\mathcal{T}) = \{ 2, 3, 4 \}$ and $\maj(\mathcal{T}) = 2 + 3 + 4 = 9$.
\end{example}

\begin{prop}\label{prop:q-ps}
    For any $n \ge 2$ and $D \subseteq [n-1]$, the specialization of $\Qtil_{n,D}$ given by setting $x_i = q^{i-1}, y_i = tq^{i-1}$ is
    \[
    \Qtil_{n,D}(1,q,q^2,\ldots;t,tq,tq^2,\ldots) = \frac{1}{(q;q)_n} \sum_{S \subseteq [n]} q^{\comaj(D,S)}t^{|S|},
    \]
    where $(a;q)_n = (1-a)(1-aq) \cdots (1-aq^{n-1})$ denotes the $q$-Pochhammer symbol.
\end{prop}
\begin{proof}
    Given a sequence $\mathbf{a} = (a_1 \le a_2 \le \cdots \le a_n)$ from the alphabet $\Acal$ which satisfies:
    \begin{equation}\label{a-con}
a_i = a_{i+1} \in \Acal_+ \Rightarrow i \not\in D, \quad a_i = a_{i+1} \in \Acal_- \Rightarrow i \in D,
\end{equation}
write $b_1 \le b_2 \le \cdots \le b_n$ for the image under projection to $\Acal_+$, and let $S_\mathbf{a} = \{ i : a_i \in \Acal_- \}$. This sequence contributes the term $z_{a_1}z_{a_2} \cdots z_{a_n}$ to $\Qtil_{n,D}(\xb;\yb)$, so upon specializing, this term becomes $q^{b_1 + b_2 + \cdots + b_n - n}t^{|S_\mathbf{a}|}$. 

We claim that \eqref{a-con} is equivalent to the condition that:
    \begin{equation}\label{b-con}
    b_i < b_{i+1} \text{ whenever } i \in D, i+1 \not\in S_\mathbf{a}, \text{ or } i \not\in D, i \in S_\mathbf{a}.
    \end{equation}
    
    Suppose that \eqref{a-con} holds. First assume that $i \in D$ and $i+1 \not\in S_\mathbf{a}$, so that $a_{i+1} \in \Acal_+$. Then we must have either $a_i < a_{i+1}$ or $a_i = a_{i+1} \in \Acal_+$. The latter possibility is not allowed by \eqref{a-con}, so we must have $a_i < a_{i+1}$. Thus $b_i \le a_i < a_{i+1} = b_{i+1}$, so $b_i < b_{i+1}$. Now assume that $i \not\in D$ and $i \in S_\mathbf{a}$. Then $b_i < a_i < a_{i+1}$, so $b_i < b_{i+1}$. Thus \eqref{b-con} holds.
    
    Now assume that \eqref{b-con} holds. Suppose $a_i = a_{i+1} \in \Acal_+$, so that $b_i = b_{i+1}$ and $i+1 \not\in S_\mathbf{a}$. Then by assumption, we must have $i \not\in D$, as required by \eqref{a-con}. On the other hand, if $a_i = a_{i+1} \in \Acal_-$, then again we have $b_i = b_{i+1}$ and $i \in S_\mathbf{a}$, which forces $i \in D$. Therefore \eqref{a-con} and \eqref{b-con} are equivalent.

    Now, for each $1 \le j \le n$, set 
    \[
    r_j = b_j - 1 - |\{ i < j : i \in D, i+1 \not\in S_\mathbf{a} \text{ or } i \not\in D, i \in S_\mathbf{a} \}|.
    \]
    Then $0 \le r_1 \le \cdots \le r_n$ and
    \begin{align*}
    b_1 + \cdots + b_n - n 
    &= r_1 + \cdots + r_n + \comaj(D,S_\mathbf{a}),
    \end{align*}
    so that
    \begin{align*}
    \Qtil_{n,D}(1,q,q^2,\ldots;t,tq,tq^2,\ldots) &= \sum_{\substack{\mathbf{a} = (a_1 \le \cdots \le a_n) \\ \text{s.t. } \eqref{a-con} \text{ holds}}} q^{b_1 + b_2 + \cdots + b_n - n}t^{|S_\mathbf{a}|} \\
    &= \sum_{S \subseteq [n]} t^{|S|} \sum_{\substack{\mathbf{a} = (a_1 \le \cdots \le a_n) \\ \text{s.t. } \eqref{a-con} \text{ holds and } S_\mathbf{a} = S}} q^{b_1 + b_2 + \cdots + b_n - n} \\
    &= \sum_{S \subseteq [n]} t^{|S|} q^{\comaj(D,S)} \sum_{0 \le r_1 \le \cdots \le r_n} q^{r_1+\cdots+r_n} \\
    &= \frac{1}{(q;q)_n} \sum_{S \subseteq [n]} q^{\comaj(D,S)} t^{|S|}.\qedhere
    \end{align*}
\end{proof}

The analogous specialization for the super Schur functions $\stil_\lambda(\xb;\yb)$ is given by Macdonald in terms of the \emph{hook lengths} $h(r,c)$ of the cells $(r,c) \in \lambda$. See also \cite{Kratt} for a bijective proof of the following identity.
\begin{thm}[{\cite[p.~27, Ex.~5 and p.~45, Ex.~3]{Mac}}]\label{qt-hook}
   The specialization of $\stil_\lambda(\xb;\yb)$ given by setting $x_i = q^{i-1}$, $y_i = tq^{i-1}$ is given by
   \[
   \stil_\lambda(1,q,q^2,\ldots;t,tq,tq^2,...) = \prod_{(r,c) \in \lambda} \frac{q^{r-1} + tq^{c-1}}{1 - q^{h(r,c)}}.
   \]
\end{thm}
In fact, (a renormalization of) the principal specialization of $\stil_\lambda(\xb;\yb)$ may also be written as the $q,t$-generating function of the statistics $\mathrm{(co)maj}$ and $\negg$:
\begin{thm}\label{thm:s-ps}
For $\lambda \vdash n$, we have
\begin{align*}
 \stil_\lambda(1,q,q^2,\ldots;t,tq,tq^2,\ldots)
   &= \frac{1}{(q;q)_n} \sum_{\mathcal{T} \in \SYT_\pm(\lambda)} q^{\comaj(T)} t^{\negg(T)} \\
   &= \frac{1}{(q;q)_n} \sum_{\mathcal{T} \in \SYT_\pm(\lambda)} q^{\maj(T)} t^{\negg(T)}.
\end{align*}
\end{thm}
\begin{proof}
    We begin by proving the first equality. By definition $\stil_\lambda(\xb;\yb) = \sum_{T \in \SYT(\lambda)} \Qtil_{n,\Des(T)}(\xb;\yb)$, so
    \begin{align*}
        \stil_\lambda(1,q,q^2,\ldots;t,tq,tq^2,\ldots) &= \sum_{T \in \SYT(\lambda)} \Qtil_{n,\Des(T)}(1,q,q^2,\ldots;t,tq,tq^2,\ldots) \\
        &= \frac{1}{(q;q)_n} \sum_{T \in \SYT(\lambda)} \sum_{S \subseteq [n]} q^{\comaj(\Des(T),S)} t^{|S|}
    \end{align*}
    by \Cref{prop:q-ps}. Note that a standard super tableau $\mathcal{T}$ is determined by a standard tableau $T = \mathcal{T}_+$ and a set of negative entries $S = \Neg(\mathcal{T})$. We then have $\comaj(\mathcal{T}) = \comaj(\Des(T),S)$, so the first equality follows.

    Now, for any $D \subseteq [n-1]$ and $S \subseteq [n]$, let 
    \[
    D^* = \{ n-i : i \in D \} \quad \text{and} \quad S^* = \{ n+1-i : i \in S \}
    \]
    denote the \emph{reverse} of $D$ and $S$, respectively. It follows from \cite[Prop.~7.19.2]{Stan} that
    \[
    |\{ T \in \SYT(\lambda) : \Des(T) = D \}| = |\{ T \in \SYT(\lambda) : \Des(T) = D^* \}|.
    \]
    Note furthermore that for any $T \in \SYT(\lambda)$, its transpose $T' \in \SYT(\lambda')$ satisfies $\Des(T') = [n-1] \backslash \Des(T)$. Therefore, there exists a bijection $\varphi : \SYT(\lambda) \to \SYT(\lambda')$ such that
    \[
    \Des(\varphi(T)) = [n-1] \backslash \Des(T)^*
    \]
    for any $T \in \SYT(\lambda)$, which first reverses the descent set of $T$ and then transposes. 
    
    Then $\varphi$ lifts to a bijection $\widetilde{\varphi} : \SYT_\pm(\lambda) \to \SYT_\pm(\lambda')$ which satisfies
    \[
    \Des(\varphi(\mathcal{T})_+) = [n-1] \backslash \Des(\mathcal{T}_+)^* \quad \text{and} \quad \Neg(\varphi(\mathcal{T})) = [n] \backslash \Neg(\mathcal{T})^*
    \]
    for any $\mathcal{T} \in \SYT_\pm(\lambda)$, given by performing $\varphi$ on $\mathcal{T}_+$ and then reverse complementing the set of negative entries. From this, we obtain
    \[
    \comaj(\Des(\varphi(\mathcal{T})_+),\Neg(\varphi(\mathcal{T}))) = \maj(\Des(\mathcal{T}_+),\Neg(\mathcal{T})) \quad \text{and} \quad \negg(\varphi(\mathcal{T})) = n - \negg(\mathcal{T}).
    \]
    It follows from \Cref{qt-hook} that
    \[
        \stil_{\lambda'}(1,q,q^2,\ldots;t,tq,tq^2,\ldots) = \prod_{(r,c) \in \lambda} \frac{q^{c-1} + tq^{r-1}}{1-q^{h(r,c)}},
    \]
    so that
    \begin{align*}
        \stil_\lambda(1,q,q^2,\ldots;t,tq,tq^2,\ldots) &= t^n\stil_{\lambda'}(1,q,q^2,\ldots;t^{-1},t^{-1}q,t^{-1}q^2,\ldots) \\
        &= t^n \sum_{T' \in \SYT(\lambda')} \Qtil_{n,\Des(T')}(1,q,q^2,\ldots;t^{-1},t^{-1}q,t^{-1}q^2,\ldots) \\
        &= \frac{t^n}{(q;q)_n} \sum_{T' \in \SYT(\lambda')} \sum_{S \subseteq [n]} q^{\comaj(\Des(T'),S)} t^{-|S|} \\
        &= \frac{t^n}{(q;q)_n} \sum_{T \in \SYT(\lambda)} \sum_{S \subseteq [n]} q^{\maj(\Des(T),S)} t^{|S|-n} \\
        &= \frac{1}{(q;q)_n} \sum_{T \in \SYT(\lambda)} \sum_{S \subseteq [n]} q^{\maj(\Des(T),S)} t^{|S|} \\
        &= \frac{1}{(q;q)_n} \sum_{\mathcal{T} \in \SYT_\pm(\lambda)} q^{\maj(\mathcal{T})}t^{\negg(\mathcal{T})}.
    \end{align*}
\end{proof}

\begin{remark}
    We are not aware of an explicit description of the map $\varphi$ in the proof of \Cref{thm:s-ps}. See \cite[\S4]{BKS} for a description of the unique $T \in \SYT(\lambda)$ with minimal and maximal major index.
\end{remark}

\Cref{thm:maj-neg-hook} now follows easily from \Cref{thm:s-ps}.
\begin{proof}[Proof (of \Cref{thm:maj-neg-hook}).]
    By \Cref{thm:s-ps},
    \begin{align*}
        \sum_{\mathcal{T} \in \SYT_\pm(\lambda)} q^{\maj(\mathcal{T})}t^{\negg(\mathcal{T})} &= (q;q)_n \stil_\lambda(1,q,q^2,\ldots;t,tq,tq^2,\ldots) \\
        &= (1-q) \cdots (1-q^n) \prod_{(r,c) \in \lambda} \frac{q^{r-1} + tq^{c-1}}{1-q^{h(r,c)}} \\
        &= \frac{1-q}{1-q} \cdots \frac{1-q^n}{1-q} \prod_{(r,c) \in \lambda} \frac{(q^{r-1} + tq^{c-1})(1-q)}{1-q^{h(r,c)}} \\
        &= [n]_q! \prod_{(r,c) \in \lambda} \frac{q^{r-1} + tq^{c-1}}{[h(r,c)]_q}.
    \end{align*}
\end{proof}

\section{Irreducible decomposition of $\Ltil_{n,m}$}\label{sec:super-KW}
We now prove \Cref{thm:super-KW}, generalizing Kr{\'a}skiewicz--Weyman's result (\Cref{thm:KW}) in the classical case. Our proof is predominantly algebraic. By manipulating the $q,t$-hook formula (\Cref{thm:maj-neg-hook}), the Schur expansion of $\Ch(\Ltil_{n,m})$ is obtained using the generalized Cauchy identity (\Cref{lem:s-p}), along with the following technical results.

Let $\Lambda_t(\xb)$ denote the ring of symmetric functions in $\xb$ with coefficients in $\CC(t)$. For any $r \ge 1$ define an operator $\Omega^r_1 : \Lambda_t(\xb)[[q]] \to \Lambda_t(\xb)$ by:
\[
\Omega^r_1f(\xb;q,t) = \frac{1}{r} \sum_{\zeta^r = 1} \zeta^{-1} f(\xb;\zeta,t),
\]
summing over all $r$-th roots of unity $\zeta$. Then $\Omega^r_1$ extracts the terms from $f(\xb;q,t)$ whose $q$-exponents are congruent to $1$ modulo $r$; that is,
\[
\Omega^r_1f(\xb;q,t) = \sum_{\ell \in \ZZ} [q^{1+\ell r}]f(\xb;q,t).
\]
We will also require the following lemma, which can be found, for instance, in \cite{Reu}.
\begin{lemma}[\cite{Reu}, Lemma~8.11]\label{lem:reu}
    For $\lambda \vdash n$, let
    \[
    \pi_\lambda(q) = \frac{(q;q)_n}{(1-q^{\lambda_1}) \cdots (1-q^{\lambda_{\ell(\lambda)}})}.
    \]
    If $d \mid n$ and $\omega$ is a primitive $d$-th root of unity, then
    \[
    \pi_\lambda(\omega) = \begin{cases}
        0 & \text{if } \lambda \neq (d^{n/d}), \\
        z_\lambda & \text{if } \lambda = (d^{n/d}).
    \end{cases}
    \]
\end{lemma}
We now proceed with the proof.
\begin{proof}[Proof (of \Cref{thm:super-KW}).]
    Let 
    \[
    f_\lambda(\xb;q,t) = \sum_{\lambda \vdash n+m} s_\lambda(\xb) \sum_{\mathcal{T} \in \SYT_\pm(\lambda)} q^{\maj(\mathcal{T})}t^{\negg(\mathcal{T})}.
    \]
    Then it suffices to show that
    \[
    \ch(\Ltil_{n,m};\xb) = [t^m]\Omega_1^{n+m}f_\lambda(\xb;q,t).
    \]
    By \Cref{qt-hook} and \Cref{lem:s-p}, we have
    \begin{align*}
    f_\lambda(\xb;q,t) &= (q;q)_n \sum_{\lambda \vdash n+m} s_\lambda(\xb) \stil_\lambda(1,q,q^2,\ldots;t,tq,tq^2,\ldots) \\
    &= (q;q)_n \sum_{\lambda \vdash n+m} \frac{p_\lambda(\xb)}{z_\lambda} \ptil_\lambda(1,q,q^2,\ldots;t,tq,tq^2,\ldots).
    \end{align*}
    Note that
    \begin{align*}
    \ptil_\lambda(1,q,q^2,\ldots;t,tq,tq^2,\ldots) &= \prod_{j=1}^{\ell(\lambda)} (p_{\lambda_j}(1,q,q^2,\ldots) + (-1)^{\lambda_j+1}p_{\lambda_j}(t,tq,tq^2,\ldots)) \\
    &= \prod_{j=1}^{\ell(\lambda)} \left( \frac{1}{1-q^{\lambda_j}} + (-1)^{\lambda_j+1}\frac{t^{\lambda_j}}{1-q^{\lambda_j}} \right),
    \end{align*}
    so that
    \[
    f_\lambda(\xb;q,t) = \sum_{\lambda \vdash n+m} \frac{p_\lambda(\xb)}{z_\lambda} \pi_\lambda(q) \sum_{J \subseteq [\ell(\lambda)]} \prod_{j \in J} (-1)^{\lambda_j+1}t^{\lambda_j}.
    \]
    Now, for any $(n+m)$-th root of unity $\zeta$, we have by \Cref{lem:reu} that $\pi_\lambda(\zeta) = 0$ unless $\lambda = (d^\frac{n+m}{d})$ and $\zeta$ is a primitive $d$-th root of unity, in which case $\pi_\lambda(\zeta) = z_\lambda$. Thus
    \begin{align*}
    \Omega_1^{n+m}f_\lambda(\xb;q,t) &= \frac{1}{n+m} \sum_{\zeta^{n+m} = 1} \zeta^{-1} f_\lambda(\xb,\zeta,t) \\
    &= \frac{1}{n+m} \sum_{\zeta^{n+m} = 1} \zeta^{-1} \sum_{\lambda \vdash n+m} \frac{p_\lambda(\xb)}{z_\lambda} \pi_\lambda(\zeta) \sum_{J \subseteq [\ell(\lambda)]} \prod_{j \in J} (-1)^{\lambda_j+1}t^{\lambda_j} \\
    &= \frac{1}{n+m} \sum_{d \mid n+m} \sum_{\substack{\zeta \text{ a primitive} \\ \text{$d$-th root of unity}}} \zeta^{-1} p_d(\xb)^\frac{n+m}{d} \sum_{J \subseteq [\frac{n+m}{d}]} (-1)^{d|J| + |J|}t^{d|J|} \\
    &= \frac{1}{n+m} \sum_{d \mid n+m} \mu(d) p_d(\xb)^\frac{n+m}{d} \sum_{J \subseteq [\frac{n+m}{d}]} (-1)^{d|J| + |J|}t^{d|J|}.
    \end{align*}
    When we extract the coefficient of $t^m$, we may restrict the first sum to all $d \mid n+m$ such that also $d \mid m$, so that $d \mid \gcd(n,m)$. Therefore
    \begin{align*}
    [t^m]\Omega_1^{n+m}f_\lambda(q,t) &= \frac{1}{n+m} \sum_{d \mid \gcd(n,m)} \mu(d) p_d(\xb)^\frac{n+m}{d} \sum_{\substack{J \subseteq [\frac{n+m}{d}], \\ |J| = \frac{m}{d}}} (-1)^{m + \frac{m}{d}} \\
    &= \frac{1}{n+m} \sum_{d \mid \gcd(n,m)} \mu(d) p_d(\xb)^\frac{n+m}{d} (-1)^{m + \frac{m}{d}} \binom{\frac{n+m}{d}}{\frac{m}{d}} \\
    &= \ch(\Ltil_{n,m};\xb),
    \end{align*}
    where the last equality follows from \Cref{thm:super-brandt}.
\end{proof}

\section{Further directions}\label{sec:addl}

\subsection{Idempotents}

Klyachko \cite{Kly} introduced a remarkable idempotent
  \[ \mathcal{K}_n = \frac{1}{n} \sum_{\sigma \in S_n} \omega^{\maj(\sigma)} \sigma \]
where $\omega$ is a primitive $n$th root of unity. He showed $V^{\otimes n} \mathcal{K}_n = \Lcal_n(V)$. See \cite[\S4]{Gar} for a combinatorial approach. Separately, Reutenauer noted
  \[ \Trm(V) = S(\Lcal(V)) = \bigoplus_{d=0}^\infty \Lcal(V)^{(d)} \]
where $\Lcal(V)^{(d)}$ is spanned by $d$th powers of elements of $\Lcal(V)$ \cite[\S2]{Reu-idemp}. Ruetenauer gave orthogonal idempotents $\rho_n^{(d)} \in \QQ S_n$ projecting $V^{\otimes n}$ onto $\Lcal(V)^{(d)}$ for $1 \leq d \leq n$. See also \cite[Thm.~7.3]{Gar} and the surrounding discussion. It would be interesting to investigate super versions of these idempotents, which we leave as an open problem.

\begin{prob}
    Find super analogues of the Klyachko idempotent $\mathcal{K}_n$ and the Reutenauer idempotents $\rho_n^{(d)}$.
\end{prob}

\subsection{Combinatorial symmetry}

A slight variation on the Kr\'askiewicz--Weyman formula shows that the multiplicity of $S^\lambda$ in $\chi^r\uparrow_{C_n}^{S_n}$ is $\#\{T \in \SYT(\lambda) : \maj(T) \equiv_n r \}$. Since the isomorphism type of $\chi^r\uparrow_{C_n}^{S_n}$ depends only on $n$ and $\gcd(r, n)$, we have the following enumerative corollary. If $\gcd(r, n) = \gcd(s, n)$, then
\begin{equation}\label{eq:sym.1}
  \#\{T \in \SYT(\lambda) : \maj(T) \equiv_n r \} = \#\{T \in \SYT(\lambda) : \maj(T) \equiv_n s \}.
\end{equation}
A combinatorial proof of \eqref{eq:sym.1} is currently unknown. See \cite{AS} for further discussion of this problem.

Similarly, we find the multiplicity of $S^\lambda$ in $(\chi^{\mathrm{cyc}} \otimes \chi^{m/\!/2+r})\uparrow_{C_{n+m}}^{S_{n+m}}$ where
  \[ m/\!/2 \coloneqq \begin{cases}m/2 & \text{if $m$ is even} \\ 0 & \text{if $m$ is odd}\end{cases} \]
is
\begin{equation}
  \#\{\mathcal{T} \in \SYT_\pm(\lambda) : \maj(\mathcal{T}) \equiv_{n+m} r, \negg(\mathcal{T}) = m\}.
\end{equation}
The isomorphism type in this case depends only on $n+m$, $\gcd(r + m/\!/2, n+m)$, and $\binom{n+m}{m}$. In particular, if $\gcd(r + m/\!/2, n+m) = \gcd(s + m/\!/2, n+m)$, then
\begin{equation}\label{eq:sym.2}
\begin{split}
  &\#\{\mathcal{T} \in \SYT_\pm(\lambda) : \maj(\mathcal{T}) \equiv_{n+m} r, \negg(\mathcal{T}) = m\} \\
  =\ &\#\{\mathcal{T} \in \SYT_\pm(\lambda) : \maj(\mathcal{T}) \equiv_{n+m} s, \negg(\mathcal{T}) = m\}.
\end{split}
\end{equation}
and, if $n, m$ are odd,
\begin{equation}\label{eq:sym.3}
\begin{split}
  &\#\{\mathcal{T} \in \SYT_\pm(\lambda) : \maj(\mathcal{T}) \equiv_{n+m} r, \negg(\mathcal{T}) = m\} \\
  =\ &\#\{\mathcal{T} \in \SYT_\pm(\lambda) : \maj(\mathcal{T}) \equiv_{n+m} r, \negg(\mathcal{T}) = n\}.
\end{split}
\end{equation}

\begin{prob}
    Find a combinatorial proof of \eqref{eq:sym.1}, \eqref{eq:sym.2}, or \eqref{eq:sym.3}.
\end{prob}

\subsection{The degree two case}

We conclude by discussing the extension of the only other known case of Thrall's problem to the super setting. Recall that
\[
\Ch(\Lcal_{(2^d)}) = h_d[e_2] = \sum_\mu s_\mu,
\]
where the sum is over all partitions of $2d$ with even column lengths. We have $\Lcal_{(2^d)} = S^d(\Lcal_2)$, so a natural extension of this case to the free Lie superalgebra would consist in determining the Schur expansion of the characters of the super higher Lie modules
\[
\Ltil_A = S^d(\Ltil_{2,0}), \Ltil_B = \bigwedge^d(\Ltil_{1,1}), \text{ and } \Ltil_C = S^d(\Ltil_{0,2}),
\]
where the matrices $A,B$, and $C$ are given by $a_{2,0} = b_{1,1} = c_{0,2} = d$ for any $d \ge 2$ and $a_{i,j} = b_{i,j} = c_{i,j} = 0$ for all other pairs $i,j$. The case of $\Ltil_A$ is classical, and the character of $\Ltil_C$ can be computed similarly by a known plethystic identity, which can also be found in \cite[Ex.~I.8.6(a)]{Mac}:
\[
\Ch(\Ltil_C) = h_d[\Ch(\Ltil_{0,2})] = h_d[h_2] = \sum_{\nu} s_\nu,
\]
where the sum is over all partition $\nu \vdash 2d$ with even parts.

Finally, the character of $\Ltil_B$ can be written as follows:
\[
\Ch(\Ltil_B) = e_d[h_2+e_2] = \sum_{k = 0}^d e_k[h_2]e_{d-k}[e_2].
\]
The Schur expansions of $e_k[h_2]$ and $e_k[e_2]$ for arbitrary $k$ are given in \cite[Ex.~I.8.6(c,d)]{Mac}, so that the Schur expansion of $\Ch(\Ltil_B)$ can be determined via the Littlewood--Richardson rule. However, to our knowledge there is no explicit combinatorial description of the coefficients in this expansion.
%
%
%

\bibliographystyle{alpha}
\bibliography{super-thrall}

\end{document}